\documentclass[a4paper,11pt]{article}
\usepackage[utf8]{inputenc}
\usepackage[T1]{fontenc}

\usepackage{amsthm,amsmath}
\usepackage{tikz}
\usepackage{mathrsfs,amssymb,amsfonts} 
\usepackage{enumitem}
\usepackage{fullpage}
\usepackage{hyperref, enumerate}
\usepackage[babel]{microtype}
\usepackage[english]{babel}
\usepackage[capitalise]{cleveref}

\usepackage{thmtools}
\usepackage{mathtools, comment}
\usepackage{amssymb}
\usepackage[nomath]{lmodern}
\usepackage{graphicx}
\usepackage{pgf,tikz,tkz-graph,subcaption}
\usetikzlibrary{arrows,shapes}
\usetikzlibrary{decorations.pathreplacing}
\usepackage{tkz-berge}
\usepackage{enumitem}
\usepackage[normalem]{ulem}
\usepackage{hyperref}
\hypersetup{colorlinks = true, linkcolor = blue, citecolor = blue, urlcolor = blue}

\allowdisplaybreaks

\usepackage[margin=1in]{geometry}
\newtheorem{defi}{Definition}

\newtheorem{thm}[defi]{Theorem}

\newtheorem{remark}[defi]{Remark}

\newcommand*{\myproofname}{Proof}

\title{On Erd\H{o}s problem $\# 648$}

\author{Stijn Cambie\thanks{Department of Computer Science, KU Leuven Campus Kulak-Kortrijk, 8500 Kortrijk, Belgium. Supported by a postdoctoral fellowship by the Research Foundation Flanders (FWO) with grant number 1225224N. Email: \protect\href{mailto:stijn.cambie@hotmail.com}{\protect\nolinkurl{stijn.cambie@hotmail.com}}}}

\begin{document}
\parindent=0cm
\maketitle




In this short article, we study a problem posed by Erd\H{o}s $30$ years ago in a collection of open problems~\cite{Erdos95}.
It was a recent problem at that point, just one year before his death, making it one of Erd\H{o}s's last problems.
The problem is listed in the online database of Erd\H{o}s problems as problem $\# 648$ (as of the time of publication),~\url{https://www.erdosproblems.com/648}.
In~\cite{Erdos95}, it was asked to estimate, as well as one can, the function $g(n)$ of the longest sequence of integers bounded by $n$ for which the smaller integers have a larger largest prime divisor.

Formally, $g(n)$ is the largest $t$ for which there exists integers $0<a_1<a_2<a_3 < \ldots <a_t \le n$ for which $P(a_i)>P(a_{i+1})$ where $P(m)=\max\{ p \colon p \text{ prime and } p \mid m\}$ denotes the largest prime divisor of $m.$
Thus, with the function $P$ representing the largest prime factor of an integer, the question is to determine the longest sequence on which the function $P$ is decreasing.
Note that if the question were to ask for $P$ to be increasing, the answer would be trivially the number of primes bounded by $n$, $\pi(n)$.
Results on similar questions concerning sequences where Euler's totient function is monotone can be found in e.g.~\cite{PPT13,Tao24}.

We will prove that $g(n)$ is of the order $\sqrt{\frac{n}{\log n}}$, i.e., $g(n) = \Theta\left( \left( \frac{n}{\log n} \right)^{1/2} \right)$ or $g(n) \asymp \left( \frac{n}{\log n} \right)^{1/2}$.

Note that we use Landau notation: for two positive functions $f$ and $g$, we write $f=O(g)$ and $f=o(g)$ if there is a constant $M>0$ resp. for every constant $M>0$, $f(x)< Mf(x)$ holds for $x$ sufficiently large.
We write $f=\Theta(g)$ if both $f=O(g)$ and $g=O(f)$ are true.
One may also prefer $f \ll g$ (equivalently $g \gg f$) and $f \asymp  g$ to denote resp. $f=O(g)$ and $f=\Theta(g).$ 
Finally, $f \sim g$ is the same as $(1 - o(1))g<f<(1 + o(1))g$

In the proof, we establish the complementary bounds: $g(n) \ll \left( \frac{n}{\log n} \right)^{1/2}$ and $g(n) \gg \left( \frac{n}{\log n} \right)^{1/2}$.

\begin{thm}
    $g(n)=\Theta( ( n/ \log n)^{1/2} ).$
\end{thm}

\begin{proof}[Proof upper bound:]
    By definition of divisibility, $a_i$ equals $q(a_i) P(a_i)$ for an integer $q(a_i)$.
Here $q(a_i)=\frac{a_i}{P(a_i)}<\frac{a_{i+1}}{P(a_{i+1})}= q(a_{i+1})$ for every $1\le i \le t-1$ and thus those integers are distinct.

As $q(a_i)P(a_i) = a_i \leq n$, it follows that for each $i$, at least one of the following holds:
\[
q(a_i) \leq \left( \frac{2n}{\log n} \right)^{1/2}
\quad \text{or} \quad P(a_i) \leq \left( \frac{n \log n}{2} \right)^{1/2}.
\]

Since the $q(a_i)$ are distinct, there are at most $ \left( \frac{2n}{\log n} \right)^{1/2}$ cases in which the first situation occurs.

Since the $P(a_i)$ are distinct primes, by the prime number theorem, the second situation occurs at most $ \pi \left( \left( \frac{n \log n}{2} \right)^{1/2} \right) \sim \frac{  \left( \frac{n \log n}{2} \right)^{1/2}}{ \frac 12 \log n} =\sqrt 2 \left(\frac{n}{\log n}\right)^{1/2}$ many times.
Thus 
$g(n)\lesssim  2 \sqrt 2\left(\frac{n}{\log n}\right)^{1/2} $.
\end{proof}


\begin{proof}[Proof lower bound:]
We construct an admissible sequence of $a_i$
’s with length $(2-o(1))\sqrt{n/ \log n}$ inductively.
Let $p_1 > p_2 > \cdots > p_r$ be the primes in the interval $\left( n^{1/2}, \left( n \log n \right)^{1/2} \right)$.
Here $r \sim \frac{\sqrt{ n \log(n)}}{ 0.5 \log( n) }  =2\sqrt{n/ \log n}.$

Let $a_1= p_1$ and thus $q_1=1$.
for every $i\ge 2$, let $q_i$ be the smallest integer for which $q_i p_i> q_{i-1} p_{i-1}$ and let $a_i=q_i p_i.$
In particular, we have $a_i=q_i p_i< q_{i-1} p_{i-1} + p_i=a_{i-1}+p_i$,
and thus by induction 
$a_i \le \sum_{j=1}^i p_j.$

By using a quantitative form of the prime number theorem~\cite[Thm.~6.9]{Dusart10}, $$\pi(x)= \frac{x}{\log x} \left( 1 +\frac1{\log x} + \frac{O(1)}{\log^2 x} \right)$$, and integration by parts, we find that
\[
\sum_{p \leq x} p = x \pi(x) - \sum_{ y \le x-1} \pi(y) =\frac{x^2}{2 \log x} + \frac{x^2}{4 (\log x)^2} + O\left( \frac{x^2}{(\log x)^3} \right).
\]


Set $x = \left( n \log n \right)^{1/2}$. The first main term is then
\[
\frac{x^2}{2 \log x}=\frac{n \log n}{\log (n \log n)} = \frac{n}{1 + \log \log n / \log n} = n \left( 1 - \frac{\log \log n}{\log n}  + O\left( \frac{(\log \log n)^2}{(\log n)^2} \right)\right).
\]
The second term is $
\frac{n \log n}{(\log (n \log n))^2} \sim \frac{n}{\log n}.$
The other error terms are smaller, so it follows that
\[
\sum_{p \leq \left( n \log n \right)^{1/2}} p = n \left( 1 - \frac{\log \log n}{\log n}  + O\left( \frac{1}{\log n} \right)\right).
\]
The latter is smaller than $n$ for large values of $n$.

Since $p_i> \sqrt n \ge \sqrt{a_i},$ $P(a_i)=p_i$ being the unique prime divisor larger than $\sqrt{a_i}$. Thus we have constructed a valid sequence of length $r \sim 2\sqrt{n/ \log n}.$
\end{proof}

\begin{remark}
    The precise form is likely $g(n)=(c+o(1)) \sqrt{ n / \log n}$ for a constant $2\le c \le 2 \sqrt 2$, since both bounds are easily seen to contain estimates that won't (asympotically) hold.
    The latter is similar to the follow-up question $\# 912$ from~\cite{Erdos82c}, where for a different number-theoretic quantity ($h(n)$, the number of different exponents in the prime factorisation of $n!$) the exact asymptotics are asked once the correct order was determined.
\end{remark}

Note that $c<2 \sqrt 2$ can be expected, since asymptotical sharpness of the upper bound would imply that a near-central term $b$ of the sequence would be the product of an integer $q \sim \left( \frac{2n}{\log n} \right)^{1/2}$ and a prime $ p> \left( \frac{n \log n}{2} \right)^{1/2},$ which is thus of the form $(1-o(1))n.$
Furthermore, there should be $(\sqrt 2 - o(1)) \left( \frac{n}{\log n} \right)^{1/2} $ integers in the interval $(b,n)$, whose largest primes are bounded by $\left( \frac{n \log n}{2} \right)^{1/2}$ and in decreasing order.

Since in the proof of the lower bound, $a_i-a_{i-1}$ can be very different from $p_i$, we expect $c$ to be larger than $2.$
It has been checked for large $n$ that $g(n)>2 \sqrt{ n / \log n}$ can be concluded from the given approach. The exact values are even much better, as one can expect that $P(a_i)$ do not have to be the primes in order. See~\url{https://github.com/StijnCambie/ErdosProblems/tree/main/EP648}.
Determining the exact value of $c$ (if existing) seems to be very hard.

\section*{Acknowledgement}

The author thanks the referees for careful reading, suggestions on related work and references and remarks to improve the readability of the paper.
We thank Thomas Bloom for his work with \url{https://www.erdosproblems.com} and for sharing a copy of~\cite{Erdos95}.

\end{document}